\documentclass[10pt]{amsart}
\usepackage{color}
\newtheorem{theorem}{Theorem}[section]
\newtheorem{lemma}[theorem]{Lemma}
\newtheorem{proposition}[theorem]{Proposition}
\newtheorem{corollary}[theorem]{Corollary}
\theoremstyle{definition}

\theoremstyle{remark}
\newtheorem{remark}[theorem]{Remark}

\numberwithin{equation}{section}
\begin{document}
\title{Complete hypersurfaces of constant isotropic curvature in space forms}
\author{H.A. Gururaja}
\address{Indian Institute of Science Education and Research,
	C/o Sree Rama Engineering College (Transit Campus),
	Rami Reddy Nagar, Karakambadi Road,
	Mangalam (P.O.) Tirupati - 517507.
	Andhra Pradesh, INDIA}
\email{gururaja@iisertirupati.ac.in}
%\thanks{First author is supported b
\author{Niteesh Kumar}
\address{Indian Institute of Science Education and Research,
C/o Sree Rama Engineering College (Transit Campus),
Rami Reddy Nagar, Karakambadi Road,
Mangalam (P.O.) Tirupati - 517507.
Andhra Pradesh, INDIA}
\email{niteeshkumar@students.iisertirupati.ac.in}
%\thanks{First author is supported by a research fellowship of CSIR, India.}


\begin{abstract}
We classify complete orientable hypersurfaces of constant isotropic curvature in space forms. We show that such a hypersurface has constant mean curvature only if it is an isoparametric hypersurface, and that it is minimal if and only if it is totally geodesic or it is the Clifford minimal hypersurface ${\mathbb S}^3(\frac{4c}{3})\times {\mathbb S}^1 (4c)$ in ${\mathbb S}^5 (c).$ 
\end{abstract}

%%%%%%%%%%%%%%%%%%%%%%%%%%%%%%%%%%
%%%%%%%%%%%%%%%%%%%%%%%%%%%%%%%%%%%
\maketitle
\section{Introduction}
\bigskip

Let $(M, g)$ be a Riemannian manifold of dimension $n\geq 4$ and $p\in M.$ On the complexified tangent space $T_pM\otimes\mathbb C$ we consider the complex-bilinear extension $( , )$ of $g$ and the Hermitian extension $< , >$  of $g.$ A vector $X\in T_pM\otimes\mathbb C$ is said to be isotropic if $(X, X)=0.$ A $2$-dimensional subspace $\sigma \subset T_pM\otimes\mathbb C$ is said to be totally isotropic if each $X\in \sigma$ is an isotropic vector. Let $\Re(p):(\wedge^2 T_p M, <<, >>)\rightarrow (\wedge^2 T_pM, <<, >>)$ denote the curvature operator of $M$ at $p.$ We also consider the complex-linear extension of $\Re(p)$ to $T_pM\otimes\mathbb C$ and  the Hermitian extension of $<<, >>$ to $\wedge^2 T_pM\otimes \mathbb C$ and denote them by the same symbols. We say that $(M, g)$ has constant isotropic curvature $C$ if 
$$\frac{<<\Re (p)( X\wedge Y),\  X\wedge Y>>}{<<X\wedge Y,\  X\wedge Y>>}=C$$
\noindent for every $\{X, Y\}$ spanning a $2$-dimensional isotropic  subspace $\sigma$ of $T_pM\otimes\mathbb C$ and $p\in M.$ \smallskip

In the above definition if we replace the equality by the inequality $>$ (resp. $\geq$), then we obtain the notion of positive (resp. nonnegative) isotropic curvature originally introduced by Micallef and Moore \cite{MM} which has been extensively studied in the literature (cf. \cite{B}, \cite{H}).
It has also turned out to be useful in resolving some long standing open problems in global Riemannian geometry such as the differentiable sphere theorem by Brendle and Schoen \cite{BS} using the Ricci flow technique. \smallskip

Manifolds of constant sectional curvature  have constant isotropic curvature. It is known (cf. [$6$, Th. $1.3$]) that if $n\geq 5,$ then $M$ must have constant sectional curvature if it has constant isotropic curvature. However, this is not true in dimension $4.$ A computation (see \S {$2$}) shows that the Riemannian product $\mathbb S^3(1)\times \mathbb R$ has constant isotropic curvature $2$ but it does not have constant sectional curvature (in fact, it does not even have constant Ricci curvature). 
\smallskip

In this paper we consider constant isotropic curvature manifolds from an extrinsic point of view. Let $M^{n+1}(c)$ denote the $(n+1)$-dimensional complete simply connected Riemannian manifold of constant sectional curvature $c.$ Thus,
$M^{n+1}=\mathbb S^{n+1}(c)$ if $c>0$; $M^{n+1}=\mathbb R^{n+1} $ if $c=0$ and $M^{n+1}=\mathbb H^{n+1}(c)$ if $c<0.$
Assume that $n\geq 4$ and let $M$ be a complete connected orientable embedded hypersurface of constant isotropic curvature $C$ in  $M^{n+1}(c).$ Our main aim is to give a classification of all such hypersurfaces. This question was motivated by the analogous classification for hypersurfaces of constant Ricci curvature in the Euclidean space due to T. Y. Thomas \cite{T}.\smallskip
 
Rotation hypersurfaces in space forms were introduced and investigated in detail by do Carmo and Dajczer \cite{CD} as a generalization of the classical surfaces of revolution in three-dimensional Euclidean space. The authors obtain explicit parametrizations for all the types of rotation hypersurfaces in space forms. We refer to this article for notation and complete details. Here we recall that there exist three distinct types of rotation hypersurfaces in $ {\mathbb H}^{n+1}(c)$ called spherical ($\delta=1$), parabolic ($\delta=0$) and hyperbolic ($\delta=-1$) depending on the nature of their parallels.\smallskip

Fix a real number $C>0.$ For any real number $0\leq \alpha<1,$ let $x_{\alpha}:\mathbb R\rightarrow \mathbb R$ be the positive function  given by 
$$ x_{\alpha}(s)=\sqrt{\frac{2}{C}}\sqrt{1-\alpha\sin{(\sqrt{C}s)}}.$$ 

\noindent The main results of this paper are summarized in the following theorem which gives a complete classification of constant isotropic curvature hypersurfaces in space forms.

\begin{theorem}\label{theorem}
	Assume that $n\geq 4.$ Let $M$ be a complete connected orientable embedded hypersurface in $M^{n+1}(c)$ which has constant isotropic curvature $C.$ \smallskip
	
	(i) If $n\geq 5,$ then either $M$ is a nontotally geodesic umbilical hypersurface in $M^{n+1}(c)$ or $M$ has constant sectional curvature $c.$ In the latter case, $M$ must be totally geodesic if $c>0.$\smallskip
	
	(ii) Assume that $n=4$ and $c=0.$ Then $C\geq 0.$ If $C=0,$ then either $M$ is locally isometric to $\mathbb R^4$ or it is a rotation hypersurface in $\mathbb R^5$ which is described by the profile function $x(s)=\sqrt{s^{2}+\beta}$ for some constant $\beta>0.$ If $C>0,$ then $M$ is an affine hyperplane, or a round sphere or a rotation hypersurface in ${\mathbb R}^5$ which is described by the function 
	$ x_{\alpha}$ for some $0\leq \alpha<1.$\smallskip
	
	(iii) Assume that $n=4$ and $c>0.$ Then $C>2c.$ If $2c<C<4c,$ then $M$ is a rotation hypersurface in ${\mathbb S}^5(c)$ which is described by the function $ x_{\alpha}$ for some $\alpha$ satisfying $0< \alpha<\frac{C}{2c}-1.$ If $C=4c,$ then either $M$ is totally geodesic or it is a rotation hypersurface in ${\mathbb S}^5(c)$  which is described by the function 
	$ x_{\alpha}$ for some $\alpha$ satisfying $0\leq \alpha<1.$ If $C>4c,$ then either $M$ is a nontotally geodesic umbilical hypersurface in ${\mathbb S}^5(c)$ or it is a rotation hypersurface in ${\mathbb S}^5(c)$ which is described by the function 
	$ x_{\alpha}$ for some constant $\alpha$ satisfying $0\leq \alpha<1.$\smallskip
	
	(iv) Assume that $n=4$ and $c<0.$ Then $C\geq 4c.$ If $C=4c,$ then either $M$ has constant sectional curvature $c$ or it is a rotation hypersurface in ${\mathbb H}^5(c)$ which is described by the function $$x(s)=\sqrt{\frac{2}{-C}} \sqrt{Ae^{\sqrt{-C}s}+Be^{-\sqrt{-C}s}-\delta}$$
for some non-negative constants $A$ and $B$ satisfying $A+B>\delta$ and $4AB>\delta^{2}.$ If $4c<C<0,$ then either $M$ is a nontotally geodesic umbilical hypersurface in ${\mathbb H}^5(c)$ or it is a rotation hypersurface which is descibed by the function $x$ above. If $C=0,$ then either $M$ is a nontotally geodesic umbilical hypersurface or it is a rotation hypersurface of spherical type in ${\mathbb H}^5(c)$ which is described by the function 
$$ \tilde{x}(s)=\sqrt{s^{2}+As+B}$$ 
 for some constants $A$ and $B$ satisfying $B>0$ and $\frac{A^2}{4B}<1.$ If $C>0,$ then either $M$ is a nontotally geodesic umbilical hypersurface or it is a spherical rotation hypersurface in ${\mathbb H}^5(c)$ which is described by the function $x_{\alpha}$ for some $0\leq \alpha<1.$\smallskip
 
 Conversely, all the above mentioned hypersurfaces have constant isotropic curvature $C.$
	\end{theorem}

Investigation of constant isotropic curvature hypersurfaces which have constant mean curvature leads to isoparametric hypersurfaces in the sense of E. Cartan \cite{C}.  Since the Clifford minimal hypersurface ${\mathbb S}^3(\frac{4c}{3})\times {\mathbb S}^1 (4c)$ in ${\mathbb S}^5(c)$ is the rotation hypersurface in ${\mathbb S}^5(c)$ corresponding to the profile function $x=\sqrt{\frac{3}{4c}},$ we obtain the following corollary to Theorem \ref{theorem} which in particular gives a characterization of the Clifford minimal hypersurface in ${\mathbb S}^5 (c).$ 

\begin{corollary}\label{corollary}
Let $M$ be a complete connected orientable embedded hypersurface in $M^{n+1}(c)$ which has constant isotropic curvature. If $M$ has constant mean curvature, then $M$ must be an isoparametric hypersurface. Moreover, $M$ is minimal if and only if it is totally geodesic or it is the Clifford minimal hypersurface ${\mathbb S}^3(\frac{4c}{3})\times {\mathbb S}^1 (4c)$ in ${\mathbb S}^5(c).$
	\end{corollary}

\bigskip

%%%%%%%%%%%%%%%%%%%%%%%%%%%%5
%%%%%%%%%%%%%%%%%%%%%%%%%%%%%%
\section{Preliminaries}
\bigskip

  The definition of constant isotropic curvature given in $\S{1}$ can be reformulated purely in terms of real vectors as follows. Any $2$-dimensional isotropic subspace $\sigma\subset T_pM\otimes\mathbb C$ admits a basis of the form $\{e_1+ie_2, e_3+ie_4\}$ for some orthonormal $4$-frame $\{e_1, e_2, e_3, e_4\}$ in $T_p M$ and, conversely, any such orthonormal $4$-frame gives rise to a $2$-dimensional isotropic subspace of $T_pM\otimes\mathbb C$ as above. From this it is easy to see that  $(M, g)$ has constant isotropic curvature $C$ if and only if we have 
 $$ K_{13}+K_{14}+K_{23}+K_{24}-2R_{1234}=C$$ for every orthonormal $4$-frame of tangent vectors $(e_1, e_2, e_3, e_4)$ in $M.$ Here, $K_{ij}$ denotes the sectional curvature of the $2$-plane spanned by the set $\{e_i, e_j\}$ and $R$ denotes the Riemann curvature tensor of $M.$\smallskip
 
 We now discuss some examples concerning constant isotropic curvature (CIC for short). Riemannian manifolds of constant sectional curvature furnish many examples of CIC manifolds. In dimensions strictly greater than $4$ these are the only examples of CIC manifolds. However, in dimension $4$ there are other examples of CIC manifolds as can be seen by the following proposition. \smallskip
 
\begin{proposition} 
	\noindent \textit{The Riemannian product ${\mathbb S}^3 (1)\times \mathbb R$ has CIC $2.$}
\end{proposition}
\noindent {\bf Proof.} Fix a point $(p, q)\in {\mathbb S}^3(1)\times \mathbb R$ and orthonormal bases $(e_1, e_2, e_3)$  for  $T_p {\mathbb S}^3(1)$ and $(e_4)$ for $T_q \mathbb R.$ Then $$K (e_1, e_3)+K (e_1, e_4)+K (e_2, e_3)+K (e_2, e_4)-2R (e_1, e_2, e_3, e_4)=2.$$

\noindent Here $K(e_i, e_j)$ denotes the sectional curvature of the $2$-plane spanned by $\{e_i, e_j\}$ and $R$ denotes the curvature tensor of ${\mathbb S}^3 (1)\times \mathbb R.$ Let $(f_1, f_2, f_3, f_4)$ be any orthonormal basis for $T_{(p, q)} ({\mathbb S}^3(1)\times \mathbb R)$ and assume that 
$$f_j= \sum_{i} a_{ij} e_i, \ \ 1\leq j\leq 4.$$ 

\noindent Let $A$ be the $4\times 4$ matrix $A=(a_{ij}).$ Since $A^{t}A=I$ we have $A A^{t}=I.$ 
In particular
\begin{equation}\label{2.1}
\sum_{i} {a_{4i}}^2=1.
\end{equation}
\noindent Let $R_1$ denote the curvature tensor of ${\mathbb S}^3(1).$ Then
$$R_{1} (u, v, z, w)=\langle u, z\rangle \langle v, w\rangle-\langle u, w\rangle \langle v, z\rangle$$
for any $\{u, v, z, w\} \subset T_p{\mathbb S}^3$ since ${\mathbb S}^3(1)$ has constant sectional curvature $1.$\smallskip

\noindent  Therefore
$$R(f_1, f_3, f_1, f_3)=R_1(\sum_{i=1}^{3} a_{i1}e_i, \sum_{i=1}^{3} a_{i3}e_i, \sum_{i=1}^{3} a_{i1}e_i, \sum_{i=1}^{3} a_{i3}e_i)$$

$$=(\sum_i a_{i1}^2) (\sum_i a_{i3}^2)-(\sum_i a_{i1} a_{i3})^2$$

$$=(1-a_{41}^2)(1-a_{43}^2)-a_{41}^2 a_{43}^2=1-a_{41}^2-a_{43}^2.$$\smallskip

\noindent A similar computation yields 
$$R(f_1, f_2, f_3, f_4)=0.$$\smallskip
\noindent Using equation (\ref{2.1}) we conclude that \smallskip

$K(f_1, f_3)+K(f_1, f_4)+K(f_2, f_3)+K(f_2, f_4)-2R(f_1, f_2, f_3,f_4)$\\
$=(1-a_{41}^2-a_{43}^2)+(1-a_{41}^2-a_{44}^2)+(1-a_{42}^2-a_{44}^2)+(1-a_{42}^2-a_{44}^2)$\\
$ =4-2({a_{41}}^2+a_{42}^2+a_{43}^2+a_{44}^2)=2.$

  \bigskip

 A similar computation gives the following proposition.

\begin{proposition}
\textit{The compact manifold ${\mathbb S }^3(1)\times {\mathbb S}^1$ has CIC $2.$ The product ${\mathbb S}^n (1)\times \mathbb R$ does not have CIC for any $n>4.$ The product ${\mathbb S}^2 (c)\times {\mathbb H}^2(-c)$ has CIC $0$ for any $c>0.$ Here, these manifolds are equipped with their standard product metrics.}
\end{proposition}

 Assume that $M$ is an oriented hypersurface in a space form $M^{n+1}(c).$ Let $\xi$ denote a smooth unit normal field on $M$ which is consistent with the orientation of $M$ and let $\mathcal {A}(X)=-\tilde{\nabla}_{X} \xi$ denote the shape operator of $M$ with respect to the unit normal $\xi.$ 
	If $R$ denotes the $(1, 3)$ curvature tensor of $M$ then the Gauss and Codazzi equations are given respectively by 
	\begin{equation}\label{gauss}
	R(X, Y)Z=c(X\wedge Y)Z+({\mathcal A}X\wedge {\mathcal A}Y)Z
	\end{equation}
		\noindent and 
		\begin{equation}\label{codazzi}
	(\nabla_{Y}{\mathcal A}) X=(\nabla_{X} {\mathcal A}) Y.
	\end{equation}\smallskip
	
	Further examples of CIC manifolds are provided by certain classes of rotation hypersurfaces in space forms (with their induced Riemannian metrics). These examples are constructed and discussed in more detail in \S{$3.$} Here we consider a special case which is not a Riemannian product.\smallskip
	
	\begin{proposition} 
		\noindent \textit{Let $\phi:{\mathbb R}^3\rightarrow {\mathbb R}^4$ be any orthogonal parametrization of the unit sphere in ${\mathbb R}^4$ and consider the isometric immersion $f:{\mathbb R}^4\rightarrow {\mathbb R}^5$ given by $$f(t_1, t_2, t_3, s)=(\sqrt{s^2+1} \ \phi_{1} (t_1, t_2, t_3),\cdots, \sqrt{s^2+1} \ \phi_{4}(t_1, t_2, t_3), \frac{s}{\sqrt{s^2+1}}),$$ where $\phi=(\phi_1, \phi_2, \phi_3, \phi_4).$ Then the rotation hypersurface $f({\mathbb R}^4)$ in ${\mathbb R}^5$ has CIC $0.$
		}
	\end{proposition}
	\noindent {\bf Proof.} From (cf. [$5$, Prop. $3.2$]) it follows that the rotation hypersurface $f({\mathbb R}^4)$ has principal curvatures (with respect to a suitably chosen orientation) $\lambda_1=\lambda_2=\lambda_3=\lambda=-\frac{1}{s^2+1}$ and $\lambda_4=\mu=\frac{1}{s^2+1}.$ Hence ${\lambda}^2+\lambda \mu=0.$ Let $\{e_1, e_2, e_3, e_4\}$ be an orthonormal basis which digonalizes the shape operator $\mathcal{A}$ of $f({\mathbb R}^4)$ and assume that
	$\mathcal{A}(e_i)=\lambda_i e_i,$ where $1\leq i\leq 4.$ From the Gauss equation (\ref{gauss}) (with $c=0$) it follows that 
	\begin{equation}\label{2.4}
	R(e_i, e_j, e_k, e_l)=\langle \mathcal {A}(e_i), e_l\rangle \langle \mathcal{A}(e_j), e_k\rangle-\langle \mathcal {A}(e_i), e_k\rangle \langle \mathcal{A} (e_j), e_l\rangle
	=\lambda_i \lambda_j (\delta_{i}^{l}\delta_{j}^{k}-\delta_{i}^{k}\delta_{j}^{l}),
	\end{equation}
	\noindent  where $\delta_{i}^{l}$ denotes the Kronecker symbol which is defined by setting $\delta_{i}^{l}=1$ if $i=l$ and equals $0$ otherwise. Hence
	$$K (e_1, e_3)+K (e_1, e_4)+K (e_2, e_3)+K (e_2, e_4)-2R (e_1, e_2, e_3, e_4)=0.$$\smallskip
	
	\noindent Let $\{f_1, f_2, f_3, f_4\}$ be any orthonormal basis for the tangent space and assume that 
	$$f_j= \sum_{i} a_{ij} e_i, \ \ 1\leq j\leq 4.$$  Let $A$ be the $4\times 4$ orthogonal matrix given by $A=(a_{ij}).$ From (\ref{2.4}) we obtain\medskip
	
	$k(f_1, f_3)$\\
	$=R(\sum_{i=1}^{4} a_{i1} e_i, \sum_{j=1}^{4} a_{j3} e_j, \sum_{k=1}^{4} a_{k3} e_k, \sum_{l=1}^{4} a_{l1} e_l)$\\
$=\sum_{1\leq i\neq j\leq 4} \lambda_i \lambda_j ( a_{i1}^{2} a_{j3}^{2}-a_{i1}a_{i3}a_{j1}a_{j3})$\\
	$={\lambda}^2 \sum_{1\leq i\neq j\leq 3} ( a_{i1}^2 a_{j3}^2-a_{i1}a_{i3}a_{j1}a_{j3})+\lambda \mu \sum_{i=1}^{3} ( a_{i1}^2 a_{43}^2-a_{i1}a_{i3}a_{41}a_{43})+\lambda\mu \sum_{j=1}^{3} ( a_{41}^{2}a_{j3}^{2}-a_{41}a_{43}a_{j1}a_{j3})$\\
	$={\lambda}^2( (1-a_{41}^2)(1-a_{43}^2)-a_{41}^2 a_{43}^2)+\lambda\mu ( (1-a_{41}^2) a_{43}^2+a_{41}^2a_{43}^2)+\lambda\mu ( (1-a_{43}^2) a_{41}^2+a_{41}^2 a_{43}^2)$\\
	$={\lambda}^2(1-a_{41}^2-a_{43}^2)+\lambda\mu (a_{41}^2+a_{43}^2)$\medskip
	
\noindent where, in the last two lines of the above computation, we have used the fact that $A$ is an orthogonal matrix. A similar computation yields $R(f_1, f_2, f_3, f_4)=0.$ Using the orthogonality of the matrix $A$ once more, we conclude that \smallskip

$K(f_1, f_3)+K(f_1, f_4)+K(f_2, f_3)+K(f_2, f_4)-2R(f_1, f_2, f_3,f_4)$\\
$=2({\lambda}^2+\lambda\mu)=0.$

\bigskip

%%%%%%%%%%%%%%%%%%%%%%%%%%%%%%%
%%%%%%%%%%%%%%%%%%%%%%%%%%%%%%5
\section{Proofs of the main results}
\bigskip

In this section we give proofs of Theorem \ref{theorem} and Corollary \ref{corollary}. We begin with the following lemma.
\begin{lemma}\label{lemma}
	 {Assume that $M$ is a hypersurface in $M^{n+1}(c)$ having CIC $C.$ Then $M$ has at most two distinct principal curvatures at each point. If we denote these principal curvatures by $\lambda$ and $\mu$ where $\lambda$ has multiplicity at least $n-1,$ then $ 4c+2(\lambda^2+\lambda \mu)= C.$ }
\end{lemma}

\smallskip

\noindent {\bf Proof.}  Let $\{e_i\}$ be an (local) orthonormal basis which diagonalizes the shape operator $\mathcal{A}$ of $M$ and let $\{\lambda_i\}$ denote the corresponding principal curvatures. Since $M$ has constant isotropic curvature $C$ we have 
 \begin{equation}\label{3.1}
 K_{ik}+K_{il}+K_{jk}+K_{jl}-2R_{ijkl}=C.
\end{equation}
\noindent Interchanging the indices $i$ and $j$ in (\ref{3.1}), we obtain 
 \begin{equation}\label{3.2}
 K_{jk}+K_{jl}+K_{ik}+K_{il}+2R_{ijkl}=C.
 \end{equation}
\noindent Adding equations (\ref{3.1}) and (\ref{3.2}) we see that  
 \begin{equation}\label{3.3}
 K_{ik}+K_{il}+K_{jk}+K_{jl}=C. 
 \end{equation}
 \noindent Now  interchanging $j$ and $k$ in (\ref{3.3}) gives 
  \begin{equation}\label{3.4}
  K_{ij}+K_{il}+K_{jk}+K_{kl}=C.
  \end{equation}
 \noindent  From equations (\ref{3.3}) and (\ref{3.4}) we conclude that 
 \begin{equation}\label{3.5}
K_{ik}+K_{jk}=K_{ij}+K_{kl} .
\end{equation}
 \noindent  Thus  
 \begin{equation}\label{3.6}
 K_{ik}+K_{jl}=K_{ij}+K_{kl}=K_{il}+K_{jk}.
 \end{equation}
 \noindent for any subset $\{i, j, k, l\}$ of $\{1,2,\cdots,n\}.$ On the other hand, it follows from the Gauss equation (\ref{gauss}) that 
  \begin{equation}\label{3.7}
 K_{ij}=c+\lambda_i \lambda_j.
 \end{equation}
 \noindent  Using equations (\ref{3.6}) and (\ref{3.7}) we obtain \
  \begin{equation}\label{3.8}
\lambda_i \lambda_j+\lambda_k \lambda_l=\lambda_i \lambda_k+\lambda_j \lambda_l=\lambda_i \lambda_l+\lambda_j \lambda_k.
\end{equation}
   \noindent  From equation (\ref{3.8}) it is easy to see  that $M$ has at most two distinct principal curvatures at each point on $M.$
  By rearranging the indices if necessary, we may assume that $\lambda_1=\lambda_2=\cdots =\lambda_{n-1}=\lambda$ and $\lambda_{n}=\mu.$ Equation (\ref{3.3}) then reads
   \begin{equation}\label{3.9}
 4c+2(\lambda^2+\lambda \mu)= C.
 \end{equation}

\begin{remark}
	The above lemma also follows from a well known theorem of E. Cartan on conformally flat hypersurfaces in space forms. It is known (cf. [$6,$  Th. $1.3$]) that CIC manifolds must be conformally flat.
\end{remark}

\begin{proposition} Assume that $n\geq 5$ and $M$ is a complete hypersurface in $M^{n+1}(c)$ which has CIC $C.$ Then either $M$ is a nontotally geodesic umbilical hypersurface in $M^{n+1}(c)$ or $M$ has constant sectional curvature $c.$
\end{proposition}

\noindent {\bf Proof.} Since $n\geq 5,\ M$ must have constant sectional curvature as observed in \cite{DN}. It then follows from equation (\ref{3.7}) that $\lambda^2=\lambda \mu$ must be a constant, and therefore $\lambda$ must be a constant on $M.$ Now consider the following cases.\smallskip

\noindent (i) $\lambda\neq 0.$ Then $\mu=\lambda$ and thus $M$ must be a nontotally geodesic umbilical hypersurface.
Proposition 1.20 in \cite{DT} (for e.g.) gives a complete classification of such hypersurfaces. In particular, if $c=0$ then $M$ must be a round  sphere in $\mathbb R^{n+1}.$

\noindent (ii)  $\lambda =0.$ Then $M$ has constant sectional curvature $c.$ In particular, if $c=0$ then  $M$ must be flat and, by a result of Hartman and Nirenberg, it must be a complete cylinder over a curve. It is also known that $M$ must be totally geodesic if $c>0.$ The case $c<0$ is considered in \cite{AP}.\smallskip

Conversely, since $M$ has constant sectional curvature in both of the above cases, $M$ must have CIC $C.$
\bigskip

We now proceed with the classification of $4$-dimensional complete CIC hypersurfaces in space forms. For this, fix a real number $C>0.$ For any real number $0\leq \alpha<1,$ let $x_{\alpha}:\mathbb R\rightarrow \mathbb R$ be the positive function  given by 
$$ x_{\alpha}(s)=\sqrt{\frac{2}{C}}\sqrt{1-\alpha\sin{(\sqrt{C}s)}}.$$
\medskip

%%%%%%%%%%%%%%%%%5
%%%%%%%%%%%%%%%%%%
\begin{subsection}
\textit{\bf  Complete CIC hypersurfaces in $\mathbb{R}^{5}.$}
\bigskip

The following proposition gives a classification of complete CIC hypersurfaces in ${\mathbb R}^5.$ \smallskip

\begin{proposition}\label{1}
Let $M$ be a complete connected orientable embedded hypersurface in ${\mathbb R}^5$ which has CIC $C.$ Then $C\geq 0.$ If $C=0,$ then either $M$ is locally isometric to $\mathbb R^4$ or it is a rotation hypersurface in $\mathbb R^5$ which is described by the profile function $x(s)=\sqrt{s^{2}+\beta}$ for some constant $\beta>0.$ If $C>0,$ then $M$ is an affine hyperplane, or a round sphere or a rotation hypersurface in ${\mathbb R}^5$ which is described by the function 
$ x_{\alpha}$ for some $0\leq \alpha<1.$ 
\end{proposition}
\noindent{\bf Proof.} Assume that $M$ is a CIC hypersurface in ${\mathbb R}^5.$ From (\ref{3.9}) we have
$\lambda^{2}+\lambda\mu=\frac{C}{2}.$ We claim that $C\geq 0.$ Otherwise $\lambda^{2}+\lambda\mu=\frac{C}{2} < 0.$
Thus $\lambda\neq 0,\  \lambda\neq\mu$ on $M$ and $\mu$ is a differentiable function of $\lambda.$
 Therefore (cf. [$5$, Th. $4.2$]) it follows that $M$ must be a rotation hypersurface in $\mathbb R^5.$
Assuming that $M$ is generated by the unit-speed profile curve described by the positive function $x=x(s):\mathbb R\rightarrow \mathbb R$ we have (cf. [$5$, Prop. $3.2$]) 
 $$ \lambda = - \frac{\sqrt{1-x'^{2}}}{x}$$ and 
 $$\mu=\frac{x''}{\sqrt{1-x'^{2}}}.$$ \smallskip
       
\noindent Substituting these values of $\lambda$ and $\mu$ in the above equation and rearranging terms we obtain the ODE
$$ (xx')'=1-\frac{C}{2} x^{2}.$$

\noindent Solving this ODE gives $$x(s)=\sqrt{\frac{2}{-C}}\sqrt{Ae^{\sqrt{-C}s}+Be^{-\sqrt{-C}s}-1} $$ \smallskip for some positive constants $A$ and $B$ satisfying $A+B>1.$ From this expression it is easy to see that the unit-speed condition $|x'|<1$ does not hold on $\mathbb R.$
\medskip

Now consider the case $C=0$ so that $\lambda^{2}+\lambda\mu=0.$ We claim that either $\lambda=0$ at each point on $M$ or $ \lambda\neq 0 $ at each point on $M.$ Suppose that $\lambda(s_{0}) \neq 0 $ for some $s_{0}\in \mathbb{R}$. By continuity there exists an open interval $I$ containg $s_0$ such that $\lambda(s)\neq 0$ for every $s\in I.$ Since $\mu\neq \lambda$ and $\mu$ is a differentiable function of $\lambda$  on $I,$ locally $M$ must be a rotation hypersurface in ${\mathbb R}^5.$ Taking $x$ to be the positive function as before we obtain the ODE 
$$\frac{1-x'^{2}}{x^{2}}-\frac{x''}{x}=0$$ which may be rewritten as $$(xx')'=1.$$

\noindent Solving this ODE (and adjusting the origin of the parameter $s$), we obtain 
$$x(s)=\sqrt{s^{2}+\beta} $$ for some constant $\beta>0.$ Moreover, $\lambda(s)=-\frac{\sqrt{\beta}}{s^{2}+\beta}$ and  $\mu(s)=\frac{\sqrt{\beta}}{s^{2}+\beta}$ in this case.\smallskip
 
\noindent Let $ D=\lbrace s\in \mathbb{R} : \lambda (s)\neq 0 \rbrace.$ Suppose $D\neq \mathbb R$ and let $D=\cup_n (a_n, b_n).$ Assume that $s_0\in (a_N, b_N).$ Since $\lim_{s\to a_{N}} \lambda(s)=-\frac{\sqrt{\beta}}{(a_{N}^{2}+\beta)}\neq 0,$ we have $\lambda(a_{N})\neq 0$ if $a_N$ is finite. Similarly $\lambda(b_{N})\neq 0$ if $b_{N}$ is finite. This shows that $D=\mathbb R$ and our claim is proved.\medskip
 
\noindent If $\lambda$ vanishes identically then $M$ must be  locally isometric to $\mathbb{R}^{4}.$ If $\lambda$ does not vanish at any point on $M$ then from the above argument it follows that $M$ must be a rotation hypersurface described by the function $x(s)=\sqrt{s^{2}+\beta}$ for some constant $\beta>0.$ 
  \medskip
 
Finally consider the case $C>0$ so that $\lambda^{2}+\lambda\mu=\frac{C}{2}>0.$ Then $\lambda\neq 0$ everywhere.
We claim that either $\lambda=\mu$ identically on $M$ or $\lambda\neq\mu$ at each point of $M.$
Assume that $\lambda (p)\neq\mu(p)$ for some $p\in\ M.$ By continuity there exists an open neighborhood of $p$ in which $\lambda$ are $\mu$ are distinct. Thus $M$ is locally a rotation hypersurface. With $x$ as before we obtain the ODE
 
 $$\frac{1-x'^{2}}{x^{2}}-\frac{x''}{x}=\frac{C}{2}$$

\noindent which, after a straighforward computation and adjustment of the origin of the parameter $s,$ is seen to have the solution
 $$ x(s)=\sqrt{\frac{2}{C}}\sqrt{1-\alpha\sin(\sqrt{C} s)}$$ for some constant $0\leq \alpha<1.$ Here the principal curvatures are given by 

 \noindent $\lambda(s)=-\frac{\sqrt{C}}{2}\sqrt{1+\frac{(1-\alpha^{2})}{(1-\alpha\sin(\sqrt{C}s))^{2}}}$ and $\mu(s)=\frac{\sqrt{C}}{2}\frac{-1+\frac{(1-\alpha^{2})}{(1-\alpha\sin(\sqrt{C}s))^{2}}}{\sqrt{1+\frac{(1-\alpha^{2})}{(1-\alpha\sin(\sqrt{C}s))^{2}}}}.$ \smallskip
 
  Let $ D=\lbrace s\in \mathbb{R} : \lambda(s)\neq \mu (s) \rbrace.$ Assume that $D\neq \mathbb R$ and let $D=\cup_n (a_n, b_n).$ Suppose that $s_0\in (a_N, b_N)$ where $x(s_0)=p.$ From the above expressions for $\lambda$ and $\mu$ we see that $\lambda(a_{N})\neq \mu (a_{N}).$  Similarly $\lambda(b_{N})\neq \mu (b_{N})$ which contradicts the fact that $D\neq \mathbb R.$ This proves the claim.\smallskip

If $M$ is umbilical then either $M$ is an affine hyperplane in $\mathbb R^5$ or a round sphere according as $\lambda=0$ or $\lambda \neq 0.$ On the other hand, if $\lambda\neq\mu$ everywhere then by the above argument it must be a rotation hypersurface in $\mathbb R^5$ which is described by the profile function 
 
 $$ x(s)=\sqrt{\frac{2}{C}}\sqrt{1-\alpha\sin{(\sqrt{C}s)}},$$ where $\alpha$ is a constant such that $0\leq \alpha<1.$\smallskip

In particular, we note from the above classification that if $\mu$ vanishes identically, then $M$ must be the Riemannian product $S^{3}({\lambda}^2)\times \mathbb R.$
\end{subsection}

%%%%%%%%%%%%%%%%%%%%%%%%%%%%%%%%%%%%%%%%%

\begin{subsection}
{Complete CIC hypersurfaces in $\mathbb{S}^{5}(c)$}\
\bigskip

The following proposition gives a classification of complete CIC hypersurfaces in ${\mathbb S}^5(c).$\smallskip

\begin{proposition}\label{2}
 Let $M$ be a complete connected orientable embedded hypersurface in ${\mathbb S}^{5}(c)$ which has CIC $C.$ Then $C>2c.$ If $2c<C<4c,$ then $M$ is a rotation hypersurface in ${\mathbb S}^5(c)$ which is described by the function $ x_{\alpha}$ for some $\alpha$ satisfying $0< \alpha<\frac{C}{2c}-1.$ If $C=4c,$ then either $M$ is totally geodesic or it is a rotation hypersurface in ${\mathbb S}^5(c)$  which is described by the function 
$ x_{\alpha}$ for some $\alpha$ satisfying $0\leq \alpha<1.$ If $C>4c,$ then either $M$ is a nontotally geodesic umbilical hypersurface in ${\mathbb S}^5(c)$ or it is a rotation hypersurface in ${\mathbb S}^5(c)$ which is described by the function 
$ x_{\alpha}$ for some constant $\alpha$ satisfying $0\leq \alpha<1.$
	\end{proposition}
	\noindent {\bf Proof.} Assume that $M$ is a complete hypersurface in ${\mathbb S}^5(c)$ of CIC $C.$ From (\ref{3.9}) we have 

\begin{equation}\label{3.10}
 {\lambda}^{2}+\lambda\mu=\frac{C-4c}{2}.
\end{equation}\smallskip

\noindent  We claim that $C> 2c$. Otherwise 
$$ \lambda^{2}+\lambda\mu=\frac{C-4c}{2}<0.$$ Thus $\lambda\neq 0,\  \lambda\neq\mu$ on $M$ and $\mu$ is a differentiable function of $\lambda.$ It follows from \cite{CD} that $M$ is a rotation hypersurface in $\mathbb S^5(c).$ Assuming that $M$ is generated by the unit-speed profile curve described by the positive function $x=x(s):\mathbb R\rightarrow \mathbb R$ we have $ \lambda = - \frac{\sqrt{1-cx^{2}-x'^{2}}}{x}$ and $\mu=\frac{x''+cx}{\sqrt{1-cx^{2}-x'^{2}}}. $ Substituting these expressions for $\lambda$ and $\mu$ in (\ref{3.10}) we obtain the ODE
$$\frac{1-cx^{2}-x'^{2}}{x^{2}}-\frac{(x''+cx)}{x}=\frac{C-4c}{2}$$
\noindent which may be rewritten as 
\begin {equation}\label{3.11}
(xx')'=1-\frac{C}{2}x^{2}.
\end{equation}\smallskip
\smallskip
       
\noindent Assuming $C<0$ and solving the ODE (\ref{3.11}) we obtain    
$$x(s)={\sqrt{\frac{2}{-C}}}\sqrt{Ae^{\sqrt{-C}s}+Be^{-\sqrt{-C}s}-1} >0$$ for some positive constants $A$ and $B$ with $A+B>1.$ From this expression we conclude as in the Euclidean case that the unit-speed condition $|x'|<1$ does not hold on $\mathbb R.$ 
\smallskip

Assuming that $C=0$ and solving the ODE we obtain $$x(s)=\sqrt{s^{2}+\beta}$$ 
for some contant $\beta>0.$ Since $cx^2+{x'}^2\geq cx^2\geq c(s^2+\beta)>1$ for all large enough values of $|s|,$ 
$\lambda(s)= - \frac{\sqrt{1-cx^{2}-x'^{2}}}{x}$ 
is not defined on $\mathbb R.$
\medskip

Now assume that $0<C\leq 2c.$ Solving the ODE (\ref{3.11}) yields the solution 
    $$x(s)={\sqrt{\frac{2}{C}}}\sqrt{1-\alpha \sin {(\sqrt{C}s)}},$$
 where $\alpha$ is a real number satisfying $0\leq \alpha<1.$ Now $cx^2+{x'}^2\geq cx^2=\frac{2c}{C}(1-\alpha \sin {(\sqrt{C}s)})=\frac{2c}{C}> 1$ at $s=0.$ Thus $\lambda$ is not defined on $\mathbb R$ and thus there do not exist complete hypersurfaces of CIC $C$ in $\mathbb{S}^{5}(c)$ if $C\leq 2c.$\medskip
       
Consider $2c<C<4c.$ In this case $$\lambda=-\sqrt{\frac{C(1-\alpha^2)} {4(1-\alpha \sin{(\sqrt{C}s)})^2}+\frac{C}{4}-c },$$ which is well-defined if and only if $$\frac{C(1-\alpha^2)}{4(1+\alpha)^2}+\frac{C}{4}-c>0$$ which, after a little algebraic manipulation, is seen to be equivalent to the inequality $\alpha<\frac{C}{2c}-1.$ Also, $\mu$ is given by
$$\mu=\frac{\frac{C(1-\alpha^{2})}{4(1-\alpha \sin{(\sqrt{C}s)})^2}-\frac{C}{4}+c}{\sqrt{\frac{C(1-\alpha^{2})}{4(1-\alpha \sin{(\sqrt{C}s)})^{2}}+\frac{C}{4}-c}}.$$\smallskip
      
In the case $C=4c$ so that $ \lambda^{2}+\lambda\mu=0$ we claim that either $\lambda=0$ at each point on $M$ or $\lambda\neq 0$ at each point on $M.$ This is proved exactly as in the Euclidean case before. If  $\lambda$ vanishes identically on $M$ then from the equation (\ref{3.7}) it follows that $M$ has constant sectional curvature $c$ and thus must be totally geodesic. On the other hand, if $\lambda$ is nowhere vanishing then the above argument shows that $M$ must be a rotation hypersurface in ${\mathbb S}^5(c)$ whose profile curve is described by the function $x=x(s): \mathbb R\rightarrow \mathbb R$ given by $$x(s)=\frac{2}{\sqrt{C}} \sqrt{1-\alpha \sin{(\sqrt{C}s)}}.$$ 
\noindent Here $0\leq\alpha <1$ is a constant.
       
 \medskip
       
   Finally assume that $C> 4c.$ Since
$\lambda^{2}+\lambda\mu=\frac{C-4c}{2}>0, \lambda\neq 0$ at each point on $M.$ Using an argument similar to the Euclidean case we can see that either $\lambda=\mu$ at each point on $M$ or $\lambda\neq\mu$ at each point on $M.$ If $\lambda=\mu$ identically, then $M$ must be a nontotally geodesic umbilical hypersurface in ${\mathbb S}^5(c),$ and as such, it must be an intersection of ${\mathbb S}^5(c)\subset {\mathbb R}^6$ with a certain affine hyperplane in ${\mathbb R}^6$ (cf. [DT, Prop. $1.20$]). If $\lambda\neq\mu$ at each point on $M,$ then $M$ must be a rotation hypersurface in ${\mathbb S}^5(c)$ which is described by the function 
$$x(s)=\sqrt{\frac{2}{C}}\sqrt{1-\alpha \sin{(\sqrt{C}s)}}$$\smallskip
 \noindent for some constant $0\leq \alpha<1.$ 
 
\end{subsection}

%%%%%%%%%%%%%%%%%%%%%%%%%%%%%%%%%
%%%%%%%%%%%%%%%%%%%%%%%%%%%%5

\begin{subsection}
{\bf Complete CIC hypersurfaces in $\mathbb{H}^{5}(c)$}\
\bigskip

The following proposition gives a classification of complete CIC hypersurfaces in ${\mathbb H}^5(c).$\smallskip

\begin{proposition}\label{3}
 Let $M$ be a complete connected orientable embedded hypersurface in ${\mathbb H}^{5}(c)$ which has CIC $C.$ Then $C\geq 4c.$ If $C=4c,$ then either $M$ has constant sectional curvature $c$ or it is a rotation hypersurface in ${\mathbb H}^5(c)$ which is described by the function $$x(s)=\sqrt{\frac{2}{-C}} \sqrt{Ae^{\sqrt{-C}s}+Be^{-\sqrt{-C}s}-\delta}$$
for some non-negative constants $A$ and $B$ satisfying $A+B>\delta$ and $4AB>\delta^{2}.$ If $4c<C<0,$ then either $M$ is a nontotally geodesic umbilical hypersurface in ${\mathbb H}^5(c)$ or it is a rotation hypersurface descibed by the function $x$ above. If $C=0,$ then either $M$ is a nontotally geodesic umbilical hypersurface or it is a rotation hypersurface of spherical type in ${\mathbb H}^5(c)$ which is described by the function 
$$ \tilde{x}(s)=\sqrt{s^{2}+As+B}$$ 
for some constants $A$ and $B$ satisfying $B>0$ and $\frac{A^2}{4B}<1.$ If $C>0,$ then either $M$ is a nontotally geodesic umbilical hypersurface or it is a spherical rotation hypersurface in ${\mathbb H}^5(c)$ which is described by the function $x_{\alpha}$ for some $0\leq \alpha<1.$
\end{proposition}
\noindent{\bf Proof.} Assume that $M$ is a hypersurface in ${\mathbb H}^{5}(c)$ of CIC $C.$
We have 
\begin{equation}\label{3.12}
\lambda^{2}+\lambda\mu=\frac{C-4c}{2}.
\end{equation}
\noindent We claim that $C\geq 4c.$ If possible, suppose that $C<4c.$ Then $\lambda^{2}+\lambda\mu=\frac{C-4c}{2}<0.$
From this it follows that $\lambda\neq 0,$\ $\lambda\neq\mu$ at each point on $M$ and that $\mu$ is a differentiable function of $\lambda.$ Thus $M$ must be a rotation hypersurface in $\mathbb{H}^{5}(c).$ With $x=x(s)$ as before, the principal curvatures $\lambda$ and $\mu$ are now given by the formulas 
$$ \lambda = - \frac{\sqrt{\delta-cx^{2}-x'^{2}}}{x}$$ and 
 
$$\mu=\frac{x''+cx}{\sqrt{\delta-cx^{2}-x'^{2}}}.$$ \smallskip
 
 \noindent Here $\delta=1,\ 0$ or $-1$ according as $M$ is spherical, parabolic or hyperbolic. We thus obtain the ODE 
$$\frac{\delta-cx^{2}-x'^{2}}{x^{2}}-\frac{(x''+cx)}{x}=\frac{C-4c}{2},$$
\noindent which may be rewritten as
\begin{equation}\label{3.13}
(xx')'=\delta-\frac{C}{2}x^{2}.
\end{equation}
\smallskip
 
\noindent Solving this ODE using the separation of variables method we obtain  $$x(s)=\sqrt{\frac{2}{-C}}\sqrt{Ae^{\sqrt{-C}s}+Be^{-\sqrt{-C}s}-\delta}$$ for some non-negative constants $A$ and $B.$  Moreover, a straightforward computation shows that $\lambda$ and $\mu$ are given by 

$$\lambda(s)=-\sqrt{-c+\frac{C}{4}-\frac{C}{4}\frac{4AB-\delta^{2}}{(Ae^{\sqrt{-C}s}+Be^{-\sqrt{-C}s}-\delta)^{2}}}$$ and

 $$\mu(s)=\frac{c-\frac{C}{4}-\frac{C}{4}\frac{4AB-\delta^{2}}{(Ae^{\sqrt{-C}s}+Be^{-\sqrt{-C}s}-\delta)^{2}}}{\sqrt{-c+\frac{C}{4}-\frac{C}{4}\frac{4AB-\delta^{2}}{(Ae^{\sqrt{-C}s}+Be^{-\sqrt{-C}s}-\delta)^{2}}}}.$$
\smallskip

\noindent From these expressions it is clear that $\lambda(s)$ and $\mu(s)$ are defined for all $s\in \mathbb{R}$ only if $-c+\frac{C}{4}\geq 0.$
\medskip

Now assume that $C=4c.$ Then $\lambda^{2}+\lambda\mu=0.$ Using an argument similar to the one that was used in the Euclidean case we see that either $\lambda$ vanishes identically or that $\lambda$ vanishes nowhere. If $\lambda$ vanishes identically then $M$ must have constant sectional curvature $c.$ If $\lambda$ is nowhere vanishing then $M$ must be a rotation hypersurface in ${\mathbb H}^5(c)$ which is described by the function
 $$x(s)=\frac{\sqrt{2}}{a}\sqrt{Ae^{as}+Be^{-as}-\delta}$$
\noindent  where the constants $A$ and $B$ satisfy $A+B>\delta$ and $4AB>\delta^{2},$ and $a=\sqrt{-C}.$ In this case, $M$ has principal curvatures

$$\lambda(s)=-\frac{\sqrt{-c(4AB-\delta^{2})}}{Ae^{as}+Be^{-as}-\delta}$$ and

 $$\mu(s)=\frac{\sqrt{-c(4AB-\delta^{2})}}{Ae^{2s}+Be^{-2s}-\delta}.$$

\smallskip

Now assume that $4c<C<0.$ Then $ \lambda^{2}+\lambda\mu=\frac{C-4c}{2}>0.$
In this case using a similar kind of argument as above we can show that either $\lambda=\mu$ identically or that $\lambda\neq\mu$ at each point on $M.$ In the first case $M$ must be a nontotally geodesic umbilical hypersurface of $\mathbb{H}^{5}(c)$ and,
in the latter case, $M$ is a rotation hypersurface in ${\mathbb H}^{5}(c)$ described by the function
 $$x(s)=\sqrt{\frac{2}{-C}} \sqrt{Ae^{\sqrt{-C}s}+Be^{-\sqrt{-C}s}-\delta} .$$
 
\noindent Moreover,

$\lambda(s)=-\sqrt{-c+\frac{C}{4}-\frac{C}{4}\frac{(4AB-\delta^{2})}{(Ae^{\sqrt{-C}s}+Be^{-\sqrt{-C}s}-\delta)^{2}}}$ and $\mu(s)=\frac{c-\frac{C}{4}-\frac{C}{4}\frac{(4AB-\delta^{2})}{(Ae^{\sqrt{-C}s}+Be^{-\sqrt{-C}s}-\delta)^{2}}}{\sqrt{-c+\frac{C}{4}-\frac{C}{4}\frac{(4AB-\delta^{2})}{(Ae^{\sqrt{-C}s}+Be^{-\sqrt{-C}s}-\delta)^{2}}}}$ 

\noindent  for all $s\in\mathbb{R}$, where $A+B>\delta$ and $4AB>\delta^{2}$. 
 \bigskip

Now consider $C=0$ so that $ \lambda^{2}+\lambda\mu=-2c.$ Then $\lambda\neq 0$ at each point on $M.$
As before we see that either $\lambda=\mu$ identically or $\lambda\neq\mu$ at each point on $M.$
 If $\lambda=\mu$ everywhere then $M$ must be nontotally geodesic umbilic hypersurface and, moreover, from (\ref{3.12}) ${\lambda}^2=-c,$ and thus by the Gauss equation the sectional curvature of $M$ vanishes identically. In the other case $M$ must be a rotation hypersurface in $\mathbb{H}^{5}(c)$  which is described by the function 
$$ x(s)=\sqrt{\delta s^{2}+As+B}$$ 

\noindent and has principal curvatures 

 $$\lambda(s)=-\sqrt{-c+\frac{4\delta B-A^{2}}{4(\delta s^{2}+As+B)}}$$ and 
 
 $$\mu(s)=\frac{c+\frac{4\delta B-A^{2}}{4(\delta s^{2}+As+B)}}{\sqrt{-c+\frac{4\delta B-A^{2}}{4(\delta s^{2}+As+B)}}}$$ 

\noindent for some constants $A$ and $B.$ For $\delta =0$ or $\delta=-1,$ the function $x$ is not defined on $\mathbb R,$ and thus there are no parabolic or hyperbolic examples in these cases. We obtain spherical examples provided $B>0$ and $\frac{A^2}{4B}<1.$

\medskip

Finally assume that $C>0.$ Since  $ \lambda^{2}+\lambda\mu=\frac{C-4c}{2}>0,$ we see that $\lambda\neq 0$ at each point on $M.$
As before we conclude that either $\lambda=\mu$ identically or $\lambda\neq\mu$ at each point on $M.$ In the first case $M$ must be a nontotally geodesic umbilical hypersurface of ${\mathbb H}^5(c).$ If $\lambda\neq\mu$ at each point on $M$ then $M$ must be a rotation hypersurface described by the function
$$x(s)=\sqrt{\frac{2}{C}}\sqrt{\delta-\alpha\sin(\sqrt{C}s)}.$$

\noindent In this case, the principal curvatures are given by  
$$\lambda(s)=-\sqrt{-c+\frac{C}{4}+\frac{C}{4}\frac{\delta^{2}-\alpha^{2}}{(\delta-\alpha \sin{(\sqrt{C}s)})^{2}}}$$
\noindent  and 
  $$\mu(s)=\frac{c-\frac{C}{4}+\frac{C}{4}\frac{\delta^{2}-\alpha^{2}}{(\delta-\alpha \sin(\sqrt{C}s)^{2}}}{\sqrt{-c+\frac{C}{4}+\frac{C}{4}\frac{\delta^{2}-\alpha^{2}}{(\delta-\alpha \sin{(\sqrt{C}s)})^2}}}.$$
 \noindent Here $0\leq\alpha<1$ is a constant.  From these expressions we conclude that $\delta =1$ and $M$ must be spherical.
\end{subsection}
\medskip

We conclude the proof of Theorem \ref{theorem} by noting that all the $4$-dimensional examples obtained in Propositions \ref{1}, \ref{2} and \ref{3}  have CIC $C$ since (cf. [$6$, Th. $1.3$]) $M$ is conformally flat and has constant scalar curvature in view of equation (\ref{3.9}).
\bigskip

\begin{subsection}
\noindent {\bf Proof of corollary \ref{corollary}}. Assume that the mean curvature $H$ of $M$ is constant. By Lemma \ref{lemma}, $H=(n-1)\lambda+\mu.$ It then follows from equation (\ref{3.9}) that $$(2-n){\lambda}^2+H\lambda-\frac{(C-4c)}{2}=0.$$ From this we conclude that $\lambda,$ and hence $\mu,$ must be constant and thus $M$ must be an isoparametric hypersurface in $M^{n+1}(c).$ If $H=0,$ then $\lambda=-\sqrt{c-\frac{C}{4}}.$ In particular, we have $C\leq 4c.$ Going through the classification in Theorem \ref{theorem} we see that $\lambda=\mu=0$ if $n\geq 5,$ or if $n=4$ and $c\leq 0.$ We also see that if $n=4$ and $c>0,$ then $M$ is either totally geodesic or it must be the rotation hypersurface in ${\mathbb S}^5(c)$ which is described by the profile function $x(s)=\sqrt{\frac{3}{4c}}.$ The latter is the Clifford minimal hypersurface ${\mathbb S}^3(\frac{4c}{3})\times {\mathbb S}^1 (4c)$ in ${\mathbb S}^5(c)$ (cf. \cite{CD}).

\end{subsection}
\bigskip
\smallskip

\maketitle
\section{Aknowledgements}
\bigskip

The authors would like to thank Abhitosh Upadhyay and Harish Seshadri for helpful comments on an earlier version of this article.

\bigskip

 %%%%%%%%%%%%%%%%%%5
 %%%%%%%%%%%%%%%%%%%%%%

\end{document}